\documentclass[12pt,a4paper,twoside]{amsart}
\usepackage{amsmath}
\usepackage{amssymb}
\usepackage{amsfonts}

\newtheorem{theorem}{Theorem}[section]
\newtheorem{lemma}[theorem]{Lemma}

\newtheorem{corollary}[theorem]{Corollary}
\theoremstyle{definition}
\newtheorem{definition}[theorem]{Definition}
\newtheorem{claim}[theorem]{Claim}
\newtheorem{example}[theorem]{Example}

\newtheorem{conjecture}[theorem]{Conjecture}
\newtheorem{remark}[theorem]{Remark}

\newcommand{\mC}{{\mathbb C}}

\newcommand{\mG}{\mathbb G}

\newcommand{\mP}{\mathbb P}

\newcommand{\mR}{{\mathbb R}}

\newcommand{\mT}{\mathbb T}

\newcommand{\ho}{\hookrightarrow}

\newcommand{\GG}{\Gamma}

\newcommand{\bo}{\omega}
\newcommand{\bl}{\lambda}
\newcommand{\bL}{\Lambda}

\newcommand{\mcC}{\mathcal C}
\newcommand{\mcD}{\mathcal D}
\newcommand{\mcE}{\mathcal E}
\newcommand{\mcF}{\mathcal F}
\newcommand{\mcG}{\mathcal G}
\newcommand{\mcH}{\mathcal H}
\newcommand{\mcI}{\mathcal I}
\newcommand{\mcJ}{\mathcal J}

\newcommand{\mcL}{\mathcal L}
\newcommand{\mcM}{\mathcal M}

\newcommand{\mcO}{\mathcal O}

\newcommand{\mcS}{\mathcal S}

\newcommand{\mcU}{\mathcal U}

\newcommand{\mcW}{\mathcal W}

\newcommand{\mcY}{\mathcal Y}
\newcommand{\mcZ}{\mathcal Z}
\newcommand{\CC}{\mathbb{C}}

\newcommand{\FF}{\mathbb{F}}

\newcommand{\WW}{\mathbb{W}}
\newcommand{\PP}{\mathbb{P}}

\newcommand{\fg}{\mathfrak g}

\newcommand{\ov}{\overline}

\newcommand{\nc}{\newcommand}
\nc{\on}{\operatorname}

\def\V{\mathcal{V}}
\def\la{\lambda}

\nc{\ol}{\overline}

\nc{\Bun}{\on{Bun}}

\nc{\mc}{\mathcal}

\newcommand{\sms}{\smallskip}
\newcommand{\ms}{\medskip}

\def\V{\mathcal{V}}

\newcommand{\of}{\overline{f}}
\newcommand{\calS}{\mathcal S}

\newcommand{\calA}{\mathcal A}

\newcommand{\calL}{\mathcal L}
\newcommand{\ome}{\omega}
\newcommand{\GL}{\operatorname{GL}}
\newcommand{\Pic}{\operatorname{Pic}}
\renewcommand{\Bun}{\operatorname{Bun}}
\newcommand{\Ad}{\operatorname{Ad}}
\newcommand{\calK}{\mathcal K}
\newcommand{\calO}{\mathcal O}
\newcommand{\Gr}{\operatorname{Gr}}
\newcommand{\Lam}{\Lambda}
\newcommand{\lam}{\lambda}
\newcommand{\oGr}{\overline{\Gr}}
\newcommand{\tGr}{\widetilde{\Gr}}
\newcommand{\grg}{\mathfrak g}
\newcommand{\calD}{\mathcal D}

\newcommand{\PGL}{\operatorname{PGL}}
\newcommand{\kap}{\kappa}
\newcommand{\Rep}{\operatorname{Rep}}
\newcommand{\calE}{\mathcal E}
\newcommand{\Hecke}{\operatorname{Hecke}}
\newcommand{\tilh}{\widetilde{h}}
\newcommand{\Aut}{\operatorname{Aut}}
\newcommand{\qlb}{\overline{\mathbb Q}_{\ell}}
\newcommand{\Fr}{\operatorname{Fr}}
\newcommand{\Tr}{\operatorname{Tr}}

\renewcommand{\Re}{\operatorname{Re}}
\newcommand{\calH}{\mathcal H}
\newcommand{\pr}{\operatorname{pr}}
\newcommand{\TT}{\mathbb T}
\newcommand\Fl{\operatorname{Fl}}
\newcommand\gam{\gamma}
\newcommand{\oc}{\overline{c}}
\newcommand{\Oper}{\operatorname{Oper}}
\newcommand{\rk}{\operatorname{rk}}
\newcommand{\DD}{\mathbb D}
\newcommand{\Spec}{\operatorname{Spec}}
\newcommand{\Sch}{\text{$\calS$ch}}
\renewcommand{\AA}{\mathbb A}
\newcommand{\SL}{\operatorname{SL}}
\newcommand\bfc{\mathbf c}
\usepackage{xcolor}
\usepackage{amscd}

\newcommand{\Hom}{\operatorname{Hom}}

\usepackage{amscd}

\begin{document}

\title[]{Automorphic functions on moduli spaces of bundles on curves over local fields: a survey}
\author{Alexander Braverman and David Kazhdan}
\address{A.B.: Department of Mathematics, University of Toronto and Perimeter Institute
of Theoretical Physics,
\newline Skolkovo Institute of Science and Technology}
\address{D.K.: Department of Mathematics, Hebrew University of Jerusalem}
\begin{abstract}
This paper is the written version of D.~Kazhdan's plenary talk at ICM 2022.
It is dedicated to an exposition of recent results and (mostly) conjectures attempting to construct an analog of the theory of automorphic functions on moduli spaces of bundles on curves over local fields (both archimedian and non-archimedian). The talk is based on joint works of D.~Kazhdan with A.Braverman, P.Etingof,  E.Frenkel and A.Polishchuk.
\end{abstract}
\maketitle
\section{Introduction}

 \subsection{Langlands correspondence over functional fields}
Let $\mcC$ be a smooth projective irreducible curve over a finite field $\FF_q$. One can consider the global field $\FF=\FF_q(X)$ of rational functions on $\mcC$ and its adele ring $\AA$. Given a split semi-simple group $G$ one can study automorphic forms on the adelic group $G(\AA)$ -- these are (by definition) irreducible representations of $G(\AA)$ which appear in the space of $\CC$-valued functions on $G(\AA)/G(\FF)$. For many purposes it is important to consider discrete automorphic representations -- these are automorphic representations appearing in $L^2(G(\AA_F/G(\FF))$.

In this introduction we restrict our attention to unramified automorphic representations -- i.e. those which have a $G( \calO_\FF)$-invariant vector where  $\calO_\FF\subset \AA$ is the ring of integral adeles. In other words we consider
functions on $G(\calO_\FF)\backslash G(\AA)/G(\FF)$ which are eigen-functions of certain commuting family of linear operators, called Hecke operators; for every place $c$ of $\FF$ (which is the same as a point of $\mcC(\overline{\FF_q})$ up to the action of Frobenius) one  constructs the algebra of Hecke operators which is isomorphic to the complexified Grothendieck ring of finite-dimensional representatons of the Langlands dual group $G^{\vee}$ (and for different $c$ these algebras commute with each other). The weak form of the Langlands conjecture (now proved by V.~Lafforgue for  global fields of positive characteristic) asserts that ( after the  replacement of  the coefficient field $\CC$ by $\qlb$,)  the common eigen-values of all the Hecke operators come from $\ell$-adic $G^{\vee}$-local systems on $\mcC$.

\sms
The quotient $G(\calO_\FF)\backslash G(\AA)/G(\FF)$ is  canonically isomorphic to the set of $\FF_q$-points of the moduli stack $\Bun_G(\mcC)$ of principal $G$-bundles on $\mcC$. Thus Hecke eigen-functions are functions on $\Bun_G(\mcC)(\FF_q)$ and unramified discrete automorphic forms correspond to Hecke eigen-function lying in $L^2(\Bun_G(X)(\FF_q))$ (with respect to the Tamagawa  measure).

\subsection{Hecke eigen-functions on moduli spaces of bundles over local fields}
This survey reports on an attempt to extend
the above constructions and results to the case when instead of a curve over $\FF_q$ we start with a curve over a local field $F$. The idea to consider Hecke eigen-functions in this case was first formulated by Langlands in the case $F=\CC$ (cf. \cite{Langl} and also \cite{Fr1}) several years ago. A systematic study of this question was started in \cite{efk1} in a slightly different framework.

\ms

Here several difficulties are present. First, since
 $\Bun$ is a stack it is not clear what space of functions on
 $\Bun(F)$
 to consider. In fact, a big part of this paper is devoted to a discussion of different vector  spaces one can attach to stacks over local fields and the relation between them (cf. Sections \ref{stacks} and \ref{bung}).
In the case when  $G$ is semi-simple and the genus $g$ of $\mcC$ is $\geq 2$ the first space is the Hilbert space $L^2(\Bun)$ which is
the completion of the space of smooth half-measures  with compact support on
$\Bun_{st}(F)$ where  $\Bun_{st}\subset \Bun$ is the open Deligne-Mumford substack
of stable bundles.

\sms

On the other hand, in \cite{efk1}
(which  deals with the case $F=\CC$) the role of Hecke operators is played by the algebra $\calD$ of global differential operators on $\Bun$ (and their complex conjugate). In fact, as was observed in \cite{bd} there is no non-trivial regular
differential operators acting on functions, but there is a  large algebra of differential operators on half-forms. This algebra $\calD$ is commutative and is equal to algebra of functions on the moduli space of certain special $G^{\vee}$-local systems on $\mcC$ called opers. This is another reason why half-forms are better suited for this problem. One of the main purposes of \cite{efk1} is a  formulation of  a conjectural  description of  eigen-values of the algebra $\calA=\calD\otimes \ov\calD$ in terms of certain $G^{\vee}$-local systems on $\mcC$ (opers with real monodromy). For $G=SL(2)$ very close conjecture was formulated by J.~Teschner in \cite{tesch}.

\ms

A systematic study of Hecke operators in the case when $F=\CC$    started in \cite{efk2}. The definition of Hecke operators in fact comes from \cite{bk} where some version of Satake isomorphism for Hecke algebras over a local field $F$ is studied (formally, \cite{bk} only deals with non-archimedian fields, but the extension to archimedian case is straightforward).
Again, it follows from \cite{bk} that in order to define Hecke operators one must work with half-forms; in this case  Hecke operators are given by certain integrals (which are not guaranteed to converge). In \cite{efk2} the authors conjectured
that these integrals in fact define
compact self-adjoint operators on $L^2(\Bun)$ for any local field $F$ (in particular, contrary to the case of finite fields, their common spectrum on $L^2(\Bun)$ is discrete); in the case $F=\CC$ it is expected that their eigen-vectors are essentially the same as the eigen-vectors for the algebra $\calA$ (we shall give a precise formulation in Section \ref{complex}). It is also explained in \cite{efk2} (in the case $F=\CC$) how to produce Hecke eigen-values from opers with real monodromy (again, this is reviewed in Section \ref{complex}). For non-archimedian fields $F$ and $G=\SL_2$ an
analogous conjectures were formulated earlier by M.Kontsevich in \cite{ko}.

\ms

 As was mentioned before, the space $L^2(\Bun)$ is not the only choice of functional space one can work with. One can define another functional space (still having to do with half-forms) on which the Hecke operators will automatically act. The relationship between this space and $L^2(\Bun)$, in the case when $G=\SL_2$, is the subject of a forthcoming paper by A.~Braverman, D.~Kazhdan and A.~Polishchuk. We review the relevant definitions and statements in Sections \ref{stacks} and \ref{bung}. As a byproduct, when $F$ is non-archimedian and the curve $\mcC$ is defined over its ring of integers $\calO_F$ (and has good reduction),  we give a conjectural construction
of an eigen-vector for the Hecke operators starting from a {\em cuspidal} eigen-function on $\Bun(\FF_q)$, where $\FF_q$ is the residue field of $F$ (but in this way one gets only a very small portion of Hecke eigen-functions).
\footnote{The construction itself is in fact not conjectural -- we can do it rigorously. But at the moment we can't prove that the resulting eigen-functions are not equal to 0.}  This is reviewed in Section \ref{hecke-local}.

We see (assuming the validity of our Conjectures) that is the non-archimedian case
at least some small part of the spectrum of Hecke operarors on $L^2(\Bun)$ has an interpretation in terms of local systems with value in the dual group. On the other had in  the archimedian case we formulate a precise conjecture of this sort in Sections \ref{complex} and \ref{real} -- in that case the whole spectrum can be (conjecturally) described in that way (following the papers \cite{efk1,efk2,efk3}).
It would be extremely interesting to relate the spectrum of Hecke operarors on $L^2(\Bun)$
to some kind of Galois data (involving the dual group $G^{\vee}$) for the non-archimedian case but at the moment we do not know how to do it
in.

\subsection{Relation of the archimedian case to geometric Langlands correspondence and conformal field theory}
In the case when the field
 $F$ is archimedian our program is related to the quantum Gauge Theory (see \cite{gw}).

\sms
In this case  Beilinson and Drinfeld associate to every $G^{\vee}$-oper $o$  a certain algebraic
$\DD$-module $M_o$ on $\Bun$ which is a {\em Hecke eigen-module} which is equipped with a canonical generator (here $\DD$ stands for the sheaf of differential operators on $\Bun$ acting on half-forms). This is an important part of a general geometric Langlands conjecture. The $\DD$-module $M_o$ can be thought of as a system of linear differential equations on $\Bun_{st}$. The corresponding Hecke eigen-half-form (in the case when $o$ has real monodromy) is a solution of both this system of equations and its complex conjugate.

\sms
The difference between the traditional categorical Langlands correspondence and the analytic Langlands correspondence for complex curves can be illustrated by an analogy with the two-dimensional conformal field theory (CFT). In CFT, there are two types of correlation functions. The first is chiral correlation functions, also known as conformal blocks. They form a vector space for fixed values of the parameters of the CFT, so we obtain a vector bundle of conformal blocks on the space of parameters, equipped with a projectively flat connection (or more generally, a twisted $\DD$-module). Conformal blocks are its {\em multi-valued} horizontal sections. The second type is the “physical” correlation functions. They can be expressed as sesquilinear combinations of conformal blocks and their complex conjugates (anti-conformal blocks), which is a {\em single-valued} function of the parameters.

\sms
The Hecke eigensheaves on $\Bun$ constructed in the categorical Langlands correspondence may be viewed as sheaves of conformal blocks of a certain CFT. They are parametrized by all $G^{\vee}$-opers on the curve. It turns out that for special $G^{\vee}$-opers (namely, the real ones) there exists a sesquilinear linear combinations of these conformal blocks and their complex conjugates which are single-valued functions (more precisely, $1/2$-measures) on $\Bun$. These are the automorphic forms of the analytic theory. Thus, the objects of the analytic theory of automorphic forms on $\Bun$ can be constructed from the objects of the categorical theory in roughly the same way as the correlation functions of CFT are obtained from conformal blocks (see \cite{Fr} and references therein for more details). An important difference with traditional CFT is that while usually in CFT the monodromy of conformal blocks is typically unitary, here the monodromy is expected to be in a split real group.

\subsection{Notations}
We shall use the letter $k$ to denote  an arbitrary field (which could be finite) and the letter $F$ for local fields. For a variety  (or stack) $X$ over $k$ we denote by $X(k)$ the set of $k$-points (for a stack we consider isomorphism classes of points).
If $F$ is non-archimedian, we denote by $\calO_F$ its ring of integers. We shall also consider the field $\calK=k((t))$ (or $F((t))$) with ring of integers which we denote just by $\calO$.

 For a split semi-simple group $G$ we denote by
  $G^{\vee}$ the Langlands dual group of $G$ considered as a group over $\CC$. We fix  a Borel subgroup  $B=TU$ of $G$, where $T$ is a maximal torus and $U$ is a maximal unipotent subgroup; similarly we have a Borel subgroup $B^{\vee} = {T^{\vee}}{U^{\vee}}\subset {G^{\vee}}$.

We denote  by  $\Lambda $ and $\Lambda ^\vee$
the lattices of coweights and of weights of $T$ (so  $\Lambda$ is also the lattice of weights of $T^{\vee}$) and by  $\bL ^+\subset \bL$ the subset of dominant coweights.

\subsection{Organization of the paper}
In Section \ref{stacks} we review some basic information about varieties and stacks over local fields and various spaces of functions on them. In Section \ref{bung} we begin the discussion of the moduli stack $\Bun$ of $G$-bundles on a curve $\mcC$ over a local field $F$ and formulate some conjectures about the relation between various function spaces one attaches to $\Bun$. In Section \ref{finite} we review the definition of Hecke operators and the formulation of the unramified Langlands correspondence for curves over $\FF_q$. In Section \ref{hecke-local} we explain the definition of Hecke operators in the case of local fields,  formulate our main conjectures and the main question specific for the non-archimedian case and also conjectures specific for the non-archimedian case.
Section \ref{complex} is dedicated to the case $F=\CC$ and Section \ref{real} to the case $F= \mR$.

\subsection{Acknowledgements}
The first author was partially supported by NSERC. The second named was partially supported by the ERC grant No 669655. We would like to thank P.Etingof, E.Frenkel and A.Polishchuk for the help with the creation of this text and D.Gaiotto and E.Witten for a very useful discussion.

\section{Smooth sections of line bundles on varieties and stacks}\label{stacks}

\subsection{Smooth sections  on varieties}

If  $X$ is  an algebraic variety over a local field $F$ (archimedian or not), the set $X(F)$ is endowed with a natural topology.
\begin{definition}
A function $f:X(F)\to \CC$ is smooth
if

a) $F$ is non-archimedian and $f$ is locally constant;

b) $F$ is archimedian and (locally) there exists a closed emebedding $X\hookrightarrow Y$ where $Y$ is a smooth variety over $F$ and a $C^{\infty}$-function $\of:Y(F)\to \CC$ such that
$f=\of|_{X(F)}$.
\end{definition}

\noindent

\sms
We denote by $C^{\infty}(X)$ the space of smooth functions on $X(F)$ and by $\calS(X)$ its subspace of functions with compact support.

\ms

For a line bundle  $\calL$ over $X$ we denote by  $\calL^0:=\calL\backslash X$ the corresponding $\mG _m$-torsor over $X$ and set
$$
|\calL|^{\kap}=\calL^0(F)\underset{F^*}\times \CC_{\kap},
$$
where $\CC_{\kap}$ denotes the 1-dimensional space $\CC$ on which $F^*$ acts by $|\cdot|^{\kap}$.
Then $|\calL|^{\kap}$ is a complex line bundle over $X(F)$. Since the bundle $ |\calL|^{\kap} $ is
 locally trivial with respect to the natural topology  we can define  its space of smooth sections which we denote by $C^{\infty}(X,|\calL|^{\kap})$. Similarly,
we denote by $\calS(X,|\calL|^{\kap})\subset C^{\infty}(X,|\calL|^{\kap}) $ the subspace of  sections with compact support.

\ms
In the case when $X$ is smooth we shall often take $\calL=\ome_X$, where $\ome_X$ is the line bundle of differential forms of top degree and write $\calS_{\kap}(X)$ instead of $\calS(X,|\ome|^{\kap})$.
The case $\kap =1/2$ is of special interest since  the space $\calS_{1/2}(X)$ is endowed with a natural Hermitian product. We denote by $L^2(X)$ its Hilbert space completion.
\begin{remark}
\begin{enumerate}
\item If $U\subset X$ is an open subset and $Z=X\backslash U$ then we have a short exact sequence
$$
0\to \calS(U,|\calL|_U^{\kap})\to \calS(X,|\calL|^{\kap})\to \calS(Z,|\calL|_Z^{\kap})\to 0.
$$
\item
More generally, instead of choosing $\calL\in \Pic(X)$ and $\kap \in \CC$ we can start with any element of $\Pic(X)\otimes \CC$ -- all the above definitions make sense in this context.
\end{enumerate}
\end{remark}

\subsection{Smooth sections on  stacks}In this
subsection we
extend the above definitions to a class of algebraic stacks.

\begin{definition}\label{adm}
 An algebraic  stack $\mcY$ is {\it admissible} if locally
there exists a  presentation of $\mcY$ as a quotient stack  $X/G $ where $X$ is a smooth variety and $G$ is an affine algebraic group. We denote by $p: X \to  \mcY$ the projection.
\end{definition}

A  presentation of $\mcY$ as a quotient $\mcY= X/\GL_n$
is called an {\it admissible presentation}.\footnote{The definition of admissibility that we use here is close to the one introduced in \cite{gk} but slightly different. It is easy to see that every admissible stack locally has an admissible presentation}

\begin{remark} \begin{enumerate}
\item Any smooth  admissible stack of finite type can be presented as a quotient  $X/\GL_n $ for a smooth variety $X$
(see \cite{gk}). As follows from the Hilbert's 90  we have $\mcY(F)= X (F) /\GL_n(F)$.
\item Any admissible stack is automatically locally of finite type.
\item A line bundle on a quotient  $X/G$  is a $G$-equivariant line bundle on $X$.

\end{enumerate}

\end{remark}
\begin{definition}
\begin{enumerate}
\item
Assume that $F$ is non-archimedian, $\mcY$ is an admissible stack of finite type over $F$. Choose an admissible presentation $\mcY=X/\GL_n$ for some variety $X$ and set
$$
\calS(\mcY,|\calL|^{\kap})=\calS(X,|\calL_X|^{\kap})_{\GL(n,F)},
$$
where the latter space stands for the space of $\GL(n,F)$-coinvariants on $\calS(X,|\calL_X|^{\kap})$.
\item
If $F$ is non-archimedian and $\mcY$ is only locally of finite type, then we can write $\mcY$ as a direct limit of open sub-stacks $\mcY_i$ of finite type over $F$
and define $\calS(\mcY,|\calL|^{\kap}):=\underset{\to}\lim\  \calS(\mcY_i,|\calL|^{\kap})$.
\item
In the case when  $F$ is archimedian we make an analogous definition but take coinvariants $\calS(X,|\calL_X|^{\kap})_{\GL(n,F)}$ in the category of topological spaces where
$\calS(X,|\calL_X|^{\kap})$ is endowed with Fr\'echet topology. \begin{footnote}{ We define the space
 $\calS(X,|\calL_X|^{\kap})_{\GL(n,F)}$ as  the quotient of $\calS(X,|\calL_X|^{\kap})$ by the closure of the subset generated by elements of
the form $g(s)-s$ where $g\in \GL_n(F)$ and $s\in \calS(X,|\calL_X|^{\kap})$.} \end{footnote}
\end{enumerate}
\end{definition}
The above definition makes sense because of the following
\begin{claim}  If  $\mcY$ is an admissible stack of finite type then the space $ \calS(\mcY,|\calL|^{\kap}) $ does not depend on a choice of an admissible  presentation $\mcY=X/\GL_n$.
\end{claim}
\begin{remark}
In the case when $F$ is non-archimedian, $\calL=\ome_X$ and $\kap =1$ this claim is proven in Section $6$ of \cite{gk}. The same arguments work in the general case.
\end{remark}
\subsection{Functoriality}
If $\mcY$ is an admissible stack and $\mcU$ is an open substack we have a natural map $\calS(\mcU,|\calL|^{\kap})\to \calS(\mcY,|\calL|^{\kap})$ which  {\color{red}TTT}
 is not be injective in general.

More generally, let $f:\mcZ\to \mcY$ be a smooth representable map of admissible stacks and $\ome_{\mcZ/\mcY}$ be the relative canonical bundle.
Then we have a natural ("integration over the fibers") map
$$
\calS(\mcZ,|\calL|^{\kap}\otimes \ome_{\mcZ/\mcY})\to \calS(\mcY,|\calL|^{\kap}).
$$

\subsection{An example: stacks over $\mcO_F$}\label{integers-eisen} In this subsection we consider the case when the field
 $F$ is non-archimedian and construct some explicit elements in $\calS(\mcY,|\calL|^{\kap})$.  Assume  that $\mcY=X/G$ where both $X$ and $G$ are defined over $\mcO_F$ and that $X_{\mcO_F}$ is a regular scheme over $\mcO_F$ such that $\mcY(F)=X(F)/G(F)$. Assume also that the line bundle $\mcL$ is defined over $\mcO_F$. Then in the same way as before we can define $\calS(\mcY_{\mcO_F},|\calL|^{\kap})$ with an obvious map $\calS(\mcY_{\mcO_F},|\calL|^{\kap})\to \calS(\mcY,|\calL|^{\kap})$.

Consider now the case when $\calL=\ome_{\mcY}$. Then the complex line bundle
 $|\calL|$ has a canonical trivialization on $\mcY(\calO_F)$.
Let $\calS(\mcY(k))$ denote the space of $\CC$-valued functions with finite support on $\mcY(k)$. Then the above trivialization gives rise to a map $\calS(\mcY(k))\to \calS(\mcY_{\mcO_F},|\calL|^{\kap})$. Composing it with the map $\calS(\mcY_{\mcO_F},|\calL|^{\kap})\to \calS(\mcY,|\calL|^{\kap})$ we get a map $E_{\mcY,\kappa}:\calS(\mcY(k))\to \calS_{\kap}(\mcY)$.
\begin{remark} \begin{enumerate}
\item
This map is often not injective.
\item We will be mostly interested in
the space $\calS_{1/2}(\mcY)$ (for a particular choice of $\mcY$).
In the case when $\mcY$ was a smooth scheme, this space had a canonical Hermitian product. We do not expect to see a Hermitian product on $\calS_{1/2}(\mcY)$ for general admissible stacks $ \mcY $
 but
we define a class of {\it excellent stacks} when such a product exists.
\item We write
 $\mcM (\mcY):=\calS_{1/2}(\mcY)$.
\end{enumerate}
\end{remark}

\subsection{Nice and  excellent stacks}
In this  subsection we assume that $\mcY$ is an admissible stack which contains
an open substack $\mcY _{vs}\subset \mcY$ such that $ \mcY _{vs} = Y_{vs}/Z$ where $Y_{vs}$ is a smooth scheme and
$Z$ is a finite group acting trivally on $Y_{vs}$.
\begin{footnote}{The subscript $vs$ stands for "very stable". The reason for this notation is that later when we work with  the stack $\Bun $ of $G$-bundles on a curve, we define  $ \Bun _{vs}\subset \Bun $ as the
open subset of very stable bundles.} \end{footnote}

\begin{remark} \begin{enumerate}
\item To simplify the notations let us assume that $Z=\{e\}$ (but generalization to arbitrary $Z$ is straightforward).
\item A choice of this open substack is not unique, and  some of the definitions below depend on this choice.
\item Let $L^2(\mcY_{vs})$ be the Hilbert space completion of the space of smooth half-measures  on $\mcY_{vs}(F)$  with compact support. It is easy to see that this space is in fact independent of the choice of $\mcY_{vs}$.
\end{enumerate}
\end{remark}

 If $\mcY$ is of finite type over $F$ we choose a presentation $\mcY=X/\GL_n$,  denote by $U$ the preimage of
 $\mcY_{vs}$ in $X$ and  by $p:U\to \mcY_{vs}$ the quotient map.

Let $s$ be a smooth section with compact support of the complex line bundle $p^*|\ome_{\mcY_{vs}}|^{\kap}\otimes |\ome_{X/\mcY}|$.
Then $s|_U$  is a section of $p^*|\ome_{\mcY_{vs}}|^{\kap}\otimes |\ome_{U/\mcY}|$. We can try to integrate it over the fibers of $p$ to get a section of $|\ome_{\mcY_{vs}}|^{\kap}$ on $\mcY_{vs}$.
The problem is that these integrals might not converge since the intersection of the support of $s$ with the fibers of the map $p$ might not be compact.

\begin{definition}\label{nice}

\begin{enumerate}
\item The stack $\mcY$ is $\kap$-bounded if there exists an open substack of finite type
 $\mcY_0\subset \mcY$ such that the map $\calS_{\kap}(\mcY_0)\to \calS_{\kap}(\mcY)$ is an isomorphism.
\item A pair $(\mcY,\mcY_{vs})$ is {\it $\kap$-nice} if $\mcY$ is $\kap$-bounded and for every $s$ as above supported on the preimage of  $\mcY _{vs}$
 the push-forward $p_*(s)$ is well-defined (i.e. it is absolutely convergent) and defines a smooth section of $|\ome_{\mcY_{vs}}|^{\kap}$ on $\mcY_{vs}$.

\item A  pair as above is {\em excellent} if it is nice for all $\kap \geq 1/2$ and for $\kap =1/2$ we have $p_*(s)\in L^2(\mcY_{vs})$ for every smooth section  $s$ with compact support.
\end{enumerate}

When the substack $\mcY _{vs}\subset \mcY$ is fixed we refer to the stack $\mcY$ as "nice" or "excellent".
\end{definition}
\begin{remark}
 The convergence in the definition of $\kappa$-niceness is automatically true for $\kappa\geq 1$.
\end{remark}

If $\mcY$ is $\kap$-nice, then the map $s\mapsto p_*(s)$ descends to a map $\calS_{\kap}(\mcY)\to C^{\infty}(\mcY_{vs})$. If $\mcY$ is excellent we get a map $\mcM(\mcY)=\calS_{1/2}(\mcY)\to L^2(\mcY_{vs})=L^2(\mcY)$.

\begin{example}Let $X =(\mP ^1)^3, G=\PGL_2$ and $\mcY =X /G$ where $G$ acts diagonally; we take $U$ to be the complement to all diagonals in $(\mP ^1)^3$. Then $G$ acts freely on $U$ and we set $\mcY_{vs}=U/G$ (note that $\mcY_{vs}$ is just Spec$F$). In this case one can check that $\mcY$ is nice for $\kap>1/3$ and
the stack $\mcY$ is excellent.
\end{example}

\section{The case of $\Bun_G$: preliminaries}\label{bung}

We  fix a split connected semi-simple group $G$ and denote by $Z$ it center.

Let $\mcC$ be  a smooth complete irreducible curve  over a field $k$.

\begin{definition}\label{Bun}
\begin{enumerate}
\item $\Bun_G$ is the stack of the  principal $G$-bundles on $\mcC$ and  $\Bun_{G,st}\subset \Bun _G$ is
the open sub-stack of stable bundles.
\item For a $G$-bundle $\mcF$ on $\mcC$ we denote by $\Ad_\mcF$ the adjoint bundle to $\mcF$ associated with the adjoint action of $G$ on $\fg$.

\item A $G$-bundle $\mcF$ is {\it very stable} if there is no non-zero  section of $ \GG (\mcC , \Ad_\mcF)\otimes \bo _\mcC $ whose values at all points of $\mcC$ are nilpotent.

\item We denote by $\Bun_{G,vs}\subset \Bun_G$ the substack of very stable bundles.
\end{enumerate}
\end{definition}
\begin{remark}If $\mcC$ is of genus  $\geq 2$ then

\begin{enumerate}
\item Every very stable bundle is stable.
\item $\Bun_{st}$ is a dense open subset of $\Bun$ of the form $Y/G$ where $Y$ is a smooth scheme of finite type over $F$
\item  $\Bun_{vs}$ is a dense open subset of $\Bun_{st}$.
\item When it does not lead to a confusion we shall drop the subscript $G$ from the notation (e.g. we shall write $\Bun$ for $\Bun_G$).
\end{enumerate}
\end{remark}

\begin{claim} \label{bounded}
The stack $  \Bun $ is $\kap$-bounded for all $\kappa$.
\end{claim}
\begin{remark}This statement is inspired by the proof of the main result of \cite{dg}.
\end{remark}

\begin{conjecture} \label{spaces}Assume that the genus $g$ of $\mcC$ is $\geq 2$.
\begin{enumerate}
\item $  \Bun $ is $\kap$-nice for $\Re(\kap)\geq 1/2$. In particular, for $\kap\geq 1/2$ we get a map $\iota_{\kappa}:
\calS_{\kappa}(\Bun)\to C^{\infty}_{\kappa}(\Bun_{vs})$.
\item For $\kappa\geq 1/2$ any section in the image of the  map $\iota_{\kappa}$ extends to a continuous section of $|\ome_{\Bun}|^{\kappa}$ on $\Bun_{st}$.

\item $  \Bun $ is excellent.
\end{enumerate}
\end{conjecture}

For $G=\PGL_2$ the 1st assertion of Conjecture \ref{spaces} (as well as some special cases of the 2nd and 3rd assertions)  will appear in a forthcoming paper of A.~Braverman, D.~Kazhdan and A.~Polishchuk. Let us note that ( again for $G=\PGL_2$) the 2nd assertion can be reduced to the following purely algebro-geometric statement using \cite{aa} (we can prove the conjecture for curves of genus 2 and 3).

\begin{conjecture}
Let $\mcE$ be a stable bundle on $\mcC$ of degree $2g-1$. Let  $F_\mcE$ denote the scheme of pairs $(\mcL,s)$ where
$\mcL\in \Pic^0(\mcC)$ and $s\in \PP (H^0(\mcC,\mcL\otimes \mcE))$.
Then
\begin{enumerate}
\item
$F_\mcE$ is  irreducible.
\item
$\dim F_{\mcE}=g.$
\item
$F_{\mcE}$ has rational singularities.
\end{enumerate}
\end{conjecture}

\section{Affine Grassmannian and Hecke operators: the case of finite field}\label{finite}
In this section we collect some facts about the canonical class of certain Schubert varieties that we shall
need in the future. All the results of this section follow easily from \cite{f} and \cite{bd}. In what follows we
a ground field $k$ and set denote by $ \calO $ the ring functions on the formal one-dimensional disc $D$ over $k$   and by $\calK$ the field of functions on the pucture disc $D^\star$. So  $\calO \sim k[[t]]$ and $\calK \sim k((t))$. We denote by $\bo _D$ the canonical bundle on $D$
and fix a square root $\bo _D^{1/2}$ (unique up to an isomorphism; the isomorphism is unique up to $\pm 1$).
%----------------------------------------------------------------------------------------------
\subsection{The affine Grassmannian}\label{grass}
Let $G$ be a split semi-simple group over $k$ and  $\Gr_G:=G(\calK)/G(\calO)$.
It is known that $\Gr_G$ has a natural structure of a proper ind-scheme over $k$ and the orbits of the group $G(\calO)$ on $\Gr_G$ are parameterised by the elements of
$\Lam_+$.

For each $\bl \in \bL ^+$ we shall denote by $\Gr^{\lam}_G$ the  corresponding orbit and by  $\oGr^{\lam}_G$ the  closure of  $\Gr^{\lam}_G$.

\subsection{Satake isomorphism} In the rest of this section we assume that $k$ is a finite field.

Let  $\mcH(G,k)$ be the algebra of compactly supported $G(\mcO)$-bi-invariant distributions on
$G(\calK)$ (by choosing a Haar measure on $G(\calK)$ such that $G(\calO)$ has volume 1, we can identify these distributions with functions). Let  $G^{\vee}$ be the Langlands dual group, considered as a group over $\CC$.
The Satake isomorphism identifies $\mcH(G,k)$ with the complexified Grothendieck ring of the category $\Rep(G^{\vee})$ of finite-dimensional representations of $G^{\vee}$. It can also be identified with the algebra $\CC[T^{\vee}]^W$ of $W$-invariant polynomial functions on $T^{\vee}$.

\subsection{Hecke operators}\label{hecke-f-subs}
Let now $\mcC$ be a smooth projective irreducible curve over $k$. As before we consider the stack $\Bun:=\Bun_G$ of principal $G$-bundles on $\mcC$. Let $c\in \mcC$ be a closed point with residue field $k'$ which is a finte extension of $k$. Choose of a local parameter near $c$,
\begin{footnote}{That is an identification of  the formal neighbourhood of $c$ with Spec$k'[[t]]$} \end{footnote}
(in the end nothing will depend on this choice)  and consider the stack $\Hecke_c$ classifying triples $(\calE_1,\calE_2,\eta)$ where every $\calE_i$ is a principal $G$-bundle on $\mcC$ and $\eta$ is an isomorphism between $\calE_1$ and $\calE_2$ on $\mcC\backslash \{ c\}$.
We have canonical projections

\begin{equation}\label{hecke-finite}
\CD
\Hecke_c @>{\pr_2}>> \Bun \\
@V{\pr_1}VV  \\
\Bun.
\endCD
\end{equation}

Every fiber of the map $\pr_2$ is isomorphic to $\Gr_G$ and this isomorphism is canonical up to the action of $G(\calO)$.
Thus every $h\in \mcH(G,k')$ defines a canonical function $\tilh$ on $\Hecke_c$. We can use it as a correspondence, and set
$$
T_{h,c}(f)=\pr_{2,*}(\pr_1^*(f)\cdot \tilh)
$$

for any $f:\Bun(\FF_q)\to \CC$.
This construction defines an action of the algebra
$\mcH(G,k')$ on the space of all functions on $\Bun(k)$ (given a choice of $c$ as above). For different choices of $c$ these operators commute.
\begin{claim}

\begin{enumerate}
\item The operators $T_{h,c}$
preserve the space $\calS(\Bun)$ of functions with finite support on $\Bun(k)$.
\item Let $L^2(\Bun(k))$ be the $L^2$-completion of the space $\calS(\Bun(k))$ with respect to the standard $L^2$-norm given by the measure on the (discrete) set $\Bun(k)$ where the volume of every $\calE$ is equal to $\frac{1}{\# \Aut(\calE)}$.
Then for every $c$ the action of $\mcH(G,k)$ extends to an action on $L^2(\Bun(k))$ by bounded operators. If $h$ is real-valued, the operator $T_{h,c}$ is self-adjoint.

\end{enumerate}
\end{claim}
\subsection{Langlands conjectures}
In the theory of automorphic forms we are usually interested in eigen-functions of all the operators $T_{h,c}$.
Let us replace the field of coefficients $\CC$ by $\qlb$ where $\ell$ is a prime number different from the characteristic of $\FF_q$. Then (the weak form) of the Langlands conjecture
states that if $f$ is such an eigen-function, then the eigen-values of all the operators $T_{h,c}$ come from
a homomorphism $\rho:\WW(\mcC)\to {G^{\vee}(\qlb)}$ where $\WW(\mcC)$ is the Weil group of $\mcC$ (a close cousin of the fundamental group of $\mcC$). In fact, in this form the Langlands conjecture has been proved by V.~Lafforgue (cf. \cite{lav}).

Let us recall the connection between Hecke-eigenvalues and homomorphisms $\rho$ as above. First of all, any $c$ defines a conjugacy class  $\Fr_c\subset  \WW(\mcC)$.
For any  $V\in \Rep(G^{\vee})$ the by Satake isomorphism associates to $V$ an element in $\mcH(G,k)$, which we denote by $h_V$. We denote $T_{V,c}$ the corresponding Hecke operator. We say that the eigen-value of an eigen-function $f$ comes from $\rho$ if
\begin{equation}\label{langlands}
T_{V,c}(f)=\Tr(\rho(\Fr_c),V)\cdot f
\end{equation}
for all $c$ and $V$.

In general, Hecke eigen-functions lie neither in $\calS(\Bun)$ nor in $L^2(\Bun)$ (here we come back to considering $\CC$-coefficients). Those which lie in the former are called {\em cuspidal}, and those which lie in the latter are called {\em discrete}. The fact that not all eigen-functions are discrete is related to the fact that the operators $T_h$ have both discrete and continuous spectrum.

\begin{remark}\label{comp}
 Note that the operators $T_{h,c}$ would be compact, if the stack $\Bun$ were of finite type over $k$ (in fact, $L^2(\Bun(k))$ would be finite-dimensional in this case), and so in that their common spectrum would be discrete. So, the
 existence of continuous spectrum of  Hecke operators is related to
the fact that $\Bun$ is not of finite type over $k$.
\end{remark}

\section{The affine Grassmannian and Hecke operators: the case of local field}\label{hecke-local}
%----------------------------------------------------------------------------------
\subsection{More on formal discs}
We are going to make a very mild change of notation (compared to the previous Section).
Namely, let $F$ be a field (very soon we shall assume that $F$ is a local field). In what follows we denote by $\calO$ some discrete valuation ring over $F$ which (as a discrete valuation ring) is isomorphic to $F((t))$ (the point is that we do not want to fix this isomorphism). We let $\calK$ be the field of fractions of $\calO$. We set $D=\Spec(\calO)$, $D^*=\Spec(\calK)$. We shall denote by $0$ the canonical $F$-point of $D$.

We let $\ome_D$ be canonical sheaf of $D$ and let $\ome_{D,0}$ be its fiber at $0$. This is a vector space over $F$.
\subsection{Line bundles on $\Gr_G$} It is well-known
(cf. \cite{bd} and \cite{f})
 that every finite-dimensional representation $V$ of
$G$ gives rise to a (determinant)
line bundle $\calL_V$ on $\Gr_G$; the fiber of this bundle over a point $g\in G(\calK)/G(\calO)$ is equal to the determinant of the vector space $g(V(\calO))/g(V(\calO))\cap V(\calO)$. In particular, we let $\calL_{\grg}$  denote the line bundle corresponding to the adjoint representation of $G$. The line bundle $ \calL_{\grg}^{-1}$ has a square
root  (unique up to isomorphism) which we denote by $\calL_{crit}$.

{\color{red} Did you define $\rho^\vee
$ ?}

The following result from \cite{bd} is crucial for us:
\begin{theorem}\label{bd}
For every $\lam\in\Lam_+$ there is a canonical isomorphism
$$
\calL_{crit}|_{\Gr_G^{\lam}}\simeq \ome_{\Gr_G^{\lam}}\otimes \ome_{D,0}^{-\langle\lambda,\rho^\vee\rangle}.
$$
(Here, as before, $\ome_{\Gr_G^{\lam}}$ denotes the canonical bundle of
on $\Gr_G^{\lam}$).
\footnote{The formulation of the theorem requires a clarification when $G$ is not simply connected, since in this case
$\langle\lambda,\rho^\vee\rangle$ might be a half-integer (and not an integer). It is sufficient for our purposes to say that we choose a square root of $\ome_{D,0}$ and that the isomorphism above is canonical up to $\pm 1$. This potential sign will disappear when we apply $|\cdot |$ to both sides which we shall do in applications.}
\end{theorem}
%---------------------------------------------------------------------

We need more information about the structure of the varieties $\oGr^{\lam}_G$.
The following result is proved in \cite{f}
(cf. also \cite{kum} and \cite{mat} for the corresponding result in characteristic 0).
\begin{theorem}\label{fal}
\begin{enumerate}
\item
Each $\oGr_G^{\lam}$ is a normal and Cohen-Macaulay projective variety over $F$.
\item
Each $\oGr_G^{\lam}$ has a resolution of singularities
\footnote{Of course, this statement is not a priori clear only if char~$F>0$.}
and for every such resolution $\pi^{\lam}:{\tGr_G^{\lam}}\to \oGr_G^{\lam}$
one has
$$
R\pi^{\lam}_*(\calO_{{\tGr_G^{\lam}}})=\calO_{\oGr_G^{\lam}}.
$$
(in other words $\oGr_G^{\lam}$ has {\it rational singularities}).
\end{enumerate}
\end{theorem}
The next result is an easy corollary of Theorem \ref{fal} and Theorem \ref{bd} (cf. \cite{bk} for a proof):
\begin{theorem}\label{canonical}
\begin{enumerate}
\item
For every $\lam\in\Lam_+$ the variety $\oGr_G^{\lam}$ is Gorenstein. Moreover,
the canonical sheaf of $\oGr_G^{\lam}$ is isomorphic to $\calL_{crit}|_{\oGr_G^{\lam}}\otimes \ome_{D,0}^{\langle\lambda,\rho^\vee\rangle}$.
Abusing the notation we shall denote this sheaf by $\ome_{\oGr_G^{\lam}}$.
\item
For any $\lam\in\Lam_+$  let $\pi^{\lam}:{\tGr_G^{\lam}}\to \oGr_G^{\lam}$ be any resolution of
singularities. Then the identification between $(\pi^{\lam})^*\ome_{\oGr_G^{\lam}}$ and
$\ome_{{\tGr_G^{\lam}}}$ that one has at the generic point of ${\tGr_G^{\lam}}$ comes from an embedding
$$
(\pi^{\lam})^*\ome_{\oGr_G^{\lam}}\hookrightarrow \ome_{{\tGr_G^{\lam}}}.
$$
(In the case $\operatorname{char}k=0$ this implies that $\oGr_G^{\lam}$ has {\it canonical singularities}).
\end{enumerate}
\end{theorem}

\subsection{Hecke algebra over local field}
In this subsection $F$ can be any local field.

Let us now work over a local field $F$ instead of $k$. So now we have $\calK=F((t)), \calO=F[[t]]$.\footnote{In the case when $F$ is non-archimedian the reader should not confuse $\calO=F[[t]]$ with $\calO_F$ which is the ring of integers of $F$.}
Then we would like to define the Hecke algebra $\mcH(G,F)$.
First we consider the space
$$
C^{\infty}_{1/2}(\Gr_G)=\underset{\leftarrow}\lim\  \calS(\oGr^{\lam}(F),|\calL_{crit}|).
$$

Assume first that $F$ is non-archimedian. Then we define $\mcH(G,F)$ to be the space of all $G(\calO)$-invariant linear functionals
on $C^{\infty}_{1/2}(\Gr_G)$ with compact support. The latter condition means that we consider functionals
$\delta: C^{\infty}_{1/2}(\Gr_G)\to \CC$ which factorize through a map $C^{\infty}_{1/2}(\Gr_G)\to \calS(\oGr^{\lam}(F),|\calL_{crit}|)$ for some $\lambda$. It is easy to see that $\mcH(G,F)$ is an algebra with respect to convolution.

It turns our that Theorems \ref{bd}, \ref{fal} and \ref{canonical} allow one to construct a lot of elements in $\calH(G,F)$ (what follows is essentially equivalent to the main result of \cite{bk}). Namely, let $\lam$ be as above and let $\phi\in C^{\infty}_{1/2}(\Gr_G)$.
Let us first trivialize the space $\ome_{D,0}$. Then $\phi|_{\Gr_G^{\lam}}$ is a distribution on $\Gr^{\lam}(F)$ and we can try to consider its integral. A priori it might not be well-defined since $\Gr_G^{\lam}(F)$ is not compact, but it is explained in \cite{bk}
that Theorems \ref{fal} and \ref{canonical} imply that in fact this integral is absolutely convergent and thus defines an element $h_{\lam}\in \mcH(G,F)$. These elements have the property that for any dominant $\lam$ and $\mu$ we have
$$
h_{\lam}\star h_{\mu}=h_{\lam+\mu}.
$$
In other words, we get an embedding $\CC[\Lam_+]\hookrightarrow \mcH(G,F)$. It is easy to see that it is actually an isomorphism.

If we do not want to trivialize the space $\ome_{D,0}$ then canonically $h_{\lambda}$ is a map from
$|\ome_{D,0}|^{-\langle\lambda,\rho^\vee\rangle}\to \mcH(G,F)$, and we get an isomorphism
$$
\bigoplus\limits_{\lambda\in \Lam_+} |\ome_{D,0}|^{\langle\lambda,\rho^\vee\rangle} \simeq \mcH(G,F)
$$
(the LHS has an obvious algebra structure).

\subsection{Hecke operators for curves over local fields: first version}
We now go back to the setup of Section \ref{bung}. We would like to define Hecke operators in this context.
First, we need to decide on what space they are going to act. The first (and the easiest) choice is to work with the space
$\mcM(\Bun)=\calS_{1/2}(\Bun)$ (another choice is discussed in the next subsection).
In what follows it will be convenient (but not necessary) to choose a particular square root $\ome^{1/2}_{\Bun}$ of $\ome_{\Bun}$ (this is always possible, but the choice is slightly not canonical).

Let us also choose a closed point $c$ of the scheme $\mcC$ with residue field $F'$ which is a finite unramified extension of $F$; we shall take $\calO$ to be the local ring of $c$ (so, it is a discrete valuation ring over $F'$ non-canonically isomorphic to $F'[[t]]$). To emphasize the dependence on $c$ we denote the corresponding Hecke algebra by $\calH_c(G)$ (instead of $\calH(G,F')$.

Then we again can consider the diagram (\ref{hecke-finite}) as in subsection \ref{hecke-f-subs}.
Then since the line bundle $\calL_{crit}$ on $\Gr_G$ is $G(\calO)$-equivariant, we can define a line bundle
$\widetilde{\calL}_{crit}$ on $\Hecke_c$ whose restriction to every fiber of $\pr_2$ is canonically isomorphic to $\calL_{crit}$ (this property makes sense since every fiber is canonically isomorphic to $\Gr_G$ up to the action of $G(\calO)$).
\begin{lemma}
We have
\begin{equation}\label{det}
\pr_1^*\ome^{1/2}_{\Bun}\simeq \pr_2^*\ome^{1/2}_{\Bun}\otimes \widetilde{\calL}_{crit}.
\end{equation}
\end{lemma}

The isomorphism (\ref{det}) easily allows one to define action of $\calH(G,F')\simeq \CC[\Lam_+]$ on $\mcM(\Bun)$.
We denote by $\TT_{\lam,c}$ the operator corresponding to $h_{\lam,c}$ (more generally, we denote by $\TT_{h,c}$ the operator corresponding to any $h\in \calH(G,F')$). For different choices of $c$ these actions commute. Therefore, one can try to study eigen-vectors of all these operators in $\mcM(\Bun)$.

\begin{remark}\label{hecke-omega}
Recall that  the operators $\TT_{\lam,c}$ are canonically defined only up to a scalar; canonically
each $\TT_{\lam,c}$ is an operator from $\mcM(\Bun)$ to $\mcM(\Bun)\otimes |\ome_{\mathcal C,c}|^{\langle\lam,\rho^{\vee}\rangle}$. Therefore when we vary $c$ each eigen-value gives rise to a section of $|\ome_{\mathcal C}|^{-\langle\lam,\rho^{\vee}\rangle}$.  This will not be important for us until the end of Section \ref{complex} (where it will in fact become quite crucial).
\end{remark}

Note that $\mcM(\Bun)$ is an analog of the space of functions with finite support on $\Bun(k)$ (where $k$ is a finite field). But unlike in the case of finite fields, we expect the following (some philosophical reasons for this difference are discussed in the next subsection):
\begin{conjecture} Assume that $F$ is non-archimedian.
Then the space $\mcM(\Bun)$ has a basis of Hecke eigen-vectors. Similarly, in the archimedian case, the space $\mcM(\Bun)$ has a topological basis of Hecke eigen-vectors.
\end{conjecture}

Before we try to say something about the eigen-values, let us discuss a slightly different version of Hecke operators.

\subsection{Hecke operators for curves over local fields: second version}
We would like to define Hecke operators in some space of actual functions on $\Bun$ (or, rather sections of $|\ome_{\Bun}|^{1/2}$), or maybe some open subset of it.
Let us assume that the genus of $\mcC$ is $\geq 2$. Then, as we have discussed before,
$\Bun$ contains a dense open sub-stack $\Bun_{st}$ of stable bundles which is a Deligne-Mumford stack.
So, one can try to start with a smooth section $\phi$ of $|\ome_{\Bun}|^{1/2}$ on $\Bun_{st}(F)$ and apply the operator $\TT_{\lam,c}$ using the diagram (\ref{hecke-finite}).
\footnote{We are slightly abusing the notation here: namely, we are going to denote by $\TT_{\lam,c}$ both the operator on $\mcM(\Bun)$ and on some space of sections of $|\ome_{\Bun}|^{1/2}$ which we are going to discuss below. We hope that it doesn't lead to a confusion}

In this case the definition will involve integration over $\Gr_G^{\lam}$ and we are not guaranteed that the corresponding integral is convergent. The trouble is caused by the following: take some $\calE\in \Bun(F)$ (which one can assume to be stable or even very stable), and consider $\pr_2^{-1}(\calE)$. Let us identify
it with $\Gr_G$ and consider the corresponding $G(\calO)$-invariant subset $\oGr^{\lam}_G$ in it. Let $S$ be a compact subset of $\Bun_{st}(F)$. Then typically $\pr_1^{-1}(S)\cap \oGr_G^{\lambda}$ is not compact.

We say that $\phi\in C^{\infty}_{1/2}(\Bun_{vs})$ is good if the integral defining $\TT_{\lam,c}(\phi)$ is absolutely convergent
and the result is again an element of $C^{\infty}_{1/2}(\Bun_{vs})$.
The following result is easy:
\begin{claim}
Assume the validity of Conjecture \ref{spaces}(1). Then the image of the map $\iota_{1/2}$ consists of
good sections and the map $\iota_{1/2}$ commutes with the operator $\TT_{\lambda,c}$.
\end{claim}
Note that the image of $\iota_{1/2}$ obviously contains $\calS_{1/2}(\Bun_{vs})$. Thus Conjecture \ref{spaces}(1) implies that any $\phi\in \calS_{1/2}(\Bun_{vs})$ is good. On the other hand, without assuming Conjecture \ref{spaces}(1) we can't a priori construct any good element of $C^{\infty}_{1/2}(\Bun_{vs})$.

We now proceed to the discussion of the action of the Hecke operators on $L^2(\Bun)$.
The main expectation is the following:

\begin{conjecture}\label{compact}
The operators $\TT_{\lam,c}$ on $L^2(\Bun)$ are bounded,
compact and self-adjoint. In particular, their common spectrum is discrete.
\end{conjecture}

Philosophically, the reason for the fact that in the case of local fields the operators $\TT_{\lam,c}$ have discrete spectrum (as opposed to the case of finite fields) is that in the case of local fields we always work only with some open subset of $\Bun$ of finite type (cf. also Claim \ref{bounded}), and as was noted in Remark \ref{comp}, the source for non-compactness of the Hecke operators in the case of finite fields has to do with the fact that the stack $\Bun$ is not globally of finite type (and in particular, not quasi-compact).

\subsection{Example}\label{mod-p}
We now want to explain how to produce some Hecke eigen-function using the construction of subsection \ref{integers-eisen}.

In the case when
$F$ is non-archimedian and that $\mcC$ is defined over $\calO_F$, i.e. we choose a model $\mcC_{\calO_F}$ of $\mcC$ over $\calO_F$. We assume that $\mcC_{\calO_F}$ is a regular scheme and we denote by $\mcC_k$ the corresponding curve over $k$. Then the stack $\Bun$ is canonically defined over $\calO_F$ and we have
the map $E_{\Bun,1/2}:\calS(\Bun(k))\to \mcM(\Bun)$( see Subsection \ref{integers-eisen}).

We claim that this map commutes with the Hecke operators in the approprite sense. Namely, let
$F'$ be a finite Galois extension of $F$ with ring of integers $\calO_{F'}$ and residue field $k'$. Then one can construct
a homomoprhism $\gam_{F'}:\calH(G,F')\to \calH(G,k')$ with the following property.
Let $c$ be a closed point of $\mcC$ whose residue field is $F'$. Note that $\mcC(F')=\mcC(\calO_{F'})$, so $c$ has canonical reduction $\oc$ which is a closed point of $\mcC_k$ with residue field $k'$.
Then for any $h\in \calH(F,F')$ and for any $\phi\in \calS(\Bun(k))$ we have
\begin{equation}\label{eis-hecke}
E_{\Bun,1/2}(T_{\gam_{F'}(h)}(\phi)=\TT_h(E_{\Bun,1/2}(\phi)).
\end{equation}
\begin{remark}
We do not know how  to describe the map $\gam_{F'}$
in general. It is easy to see  that $\gam_{F'}(h_{\lambda})$ is supported on
$\oGr_G^{\lambda}(k')$ (when viewed as a function on $\Gr_G(k')$). But this information is sufficient only in the case when $G=\PGL_n$  when minuscule coweights generate $\bL$.
\end{remark}
Equation \ref{eis-hecke} implies that $E_{\Bun,1/2}$ sends Hecke eigen-functions to Hecke eigen-functions.
This operator is certainly not injective, but we expect it to be injective on cuspidal functions.
More precisely, (  assuming the validity of Conjecture \ref{spaces})
we  formulate the following
\begin{conjecture}\label{cusp-unitary}
Assume the validity of Conjecture \ref{spaces}. Then the composition of $\iota_{1/2}\circ E_{\Bun,1/2}$ is unitary on cuspidal functions.
\end{conjecture}
Conjecture \ref{cusp-unitary} implies that we can attach a non-zero Hecke eigen-vector in $L^2(\Bun)$ to any {\em cuspidal}
Hecke eigen-function in $\calS(\Bun(k))$. On the other hand, we expect that the map $E_{\Bun,1/2}$ is highly non-injective on non-cuspidal functions. For example, let $G=\PGL_2$ and let $\calS(\Bun(k))^{\perp}_{cusp}$ denote the space of functions
with finite support which are orthogonal to all cuspidal functions (with respect to the standard Hermitian product). This space is infinite-dimensional but we expect that

$$
\dim E_{\Bun,1/2}(\calS(\Bun(k))^{\perp}_{cusp})=1.
$$

Note that equation (\ref{eis-hecke}) implies that the action of any $\TT_{\lam,c}$ on any section in the image of $E_{\Bun,1/2}$ depends only on $\oc$ (and not on $c$). This is certainly a very restrictive condition.
Also, one should think about $E_{\Bun,1/2}$ as some kind of Eisenstein series operator between the group $G(k)$ and the group $G(F)$ (with $G(\calO_F)$ playing the role of a parabolic subgroup). This is in fact the source for our notation.

%{\color{red} Shouldn't we use $p_i^\star$ instead of $\bar q_i^\star$? Is $\bar q_3=\bar q$? Is $\mX=\mcC$? What is middle $\omega$?}

\subsection{Parabolic bundles}\label{parabolic}
We would like to introduce a generalization of the above setup, which allows in particular, to consider the case of  curves of genus $\leq 1$ when we may analyze some explicit  non-trival examples.

\begin{definition}
\begin{enumerate}
\item Let us denote by $\Fl$ the variety of Borel subgroups of $G$.
\item
For a $G$-bundles $\mcF$ on $\mcC$ we denote by $\Fl_\mcF$ the associated $\Fl$-bundle over $\mcC$.
\item For a divisor $D\subset \mcC$ defined over $k$ we  denote by $\Bun^D$ the stack of $G$-bundles $\mcF$ on $\mcC$ with a section $a$ of $\Fl_\mcF$ over $D$.

\end{enumerate}

\end{definition}

It is easy to extend the definition of the Hecke operators  $\mT_{\lam,c}$  for $c\notin D$.
All our constructions and conjectures can be extended to this case. As was noted above, considering parabolic points allows one to consider
explicit examples. For example, in the case when $\mcC={\mathbb P}^1$, $D$ consists of at least 3 points and $G$ is of rank 1, Conjecture \ref{compact} is proved in \cite{efk3} (Proposition 3.13).

\subsection{Another construction of a space with Hecke action}
\subsubsection{The map $E_{\mathcal Y,\kap,n}$}
Here we would like to discuss how to generalize the construction of subsections \ref{integers-eisen} and  \ref{mod-p}. Namely, let $\mathcal Y$ be as in subsection \ref{integers-eisen}. Let $A_n=\mathcal O_F/\mathfrak m_F^n$.
Let $\mathcal Y_n$ denote the reduction of $\mathcal Y_{\mathcal O_F}$ modulo $\mathfrak m_F^n$.
This is a regular stack over $A_n$. We consider the set $\mathcal Y_n(A_n)$ and we set $\mathcal S(\mathcal Y_n(A_n))$ to be the vector space of $\mathbb C$-valued functions on $\mathcal Y_n(A_n)$ with finite support.
Then for any $\kappa\in \mathbb C$ we have the obvious map
$$
\mathcal S(\mathcal Y_n(A_n))\to \calS(\mathcal Y_{\calO_F})\simeq \calS(\mathcal Y_{\calO_F},|\ome_{\mathcal Y}^{\kappa}|).
$$
Composing it with the natural map $\calS(\mathcal Y_{\calO_F},|\ome_{\mathcal Y}^{\kappa}|)\to \calS_{\kappa}(\mathcal Y)$
we get a map
$$
E_{\mathcal Y,\kap,n}:\mathcal S(\mathcal Y_n(A_n))\to \calS_{\kappa}(\mathcal Y).
$$
\subsubsection{Commutation with Hecke operators}
We now want to specialize to the case $\mathcal Y=\Bun$ and $\kap=1/2$.
We claim that in this case the map $E_{\Bun,1/2,n}$ commutes with the Hecke operators in the sense similar to (\ref{eis-hecke}). To explain the formulation we first need to discuss an analog of the homomorphism $\gam_F$; this is a local question.

Namely, let us consider the ring $\calK_n=A_n((t))$. This is a locally compact topological ring; its subring $\calO_n=A_n[[t]]$ is open and compact. Thus the group $G(\calK_n)$ is a totally disconnected locally compact topological group with an open compact subgroup $G(\calO_n)$. Hence, we may consider the corresponding Hecke algebra
$$
\calH_n(G)=\calH(G(\calK_n),G(\calO_n)).
$$

Here is a variant. Let $\mathcal C$ be a smooth projective curve over $\calO_F$. We denote by $\mathcal C_n$ its reduction mod $\mathfrak m_F^n$. Let $F'$ be a finite unramified extension of $F$ and let $c$ be an $F'$-point of $\mathcal C$. As before we can also view it as an $\calO_{F'}$-point of $\mathcal C$ and we denote by $c_n$ its reduction modulo $\mathfrak m_{F'}^n$. This is an $A_n'$-point of $\mathcal C_n$. Then we might consider the corresponding Hecke algebra
$\calH_{c_n}(G)$. It is non-canonically isomorphic to $\calH_n(G)$.

This Hecke algebra is quite bad: it is not commutative for $n>1$ and apparently it does not have any reasonable description.
However, it has the following two important features:

1) Let $\mathcal C$ above and let $c$ be a point of $\mathcal C$ defined over a finite unramified extension $F'$ of $F$. Then the (non-commutative) algebra $\calH_{c_n}(G)$ acts on $\calS(\Bun_n)$. Given $h\in \calH_{c_n}(G)$ we denote by $T_{h,n}$ the corresponding operator on $\calS(\Bun_n)$.

2) We have a canonical homomorphism $\gam_{c,n}:\calH_{c}(G)\to \calH_{c_n}(G)$.

3) For any $h\in \calH_{c}(G)$ and $\phi\in \calS(\Bun_n)$ we have
\begin{equation}\label{eisen-n}
E_{\Bun,1/2,n}(\gam_{c,n}(\TT_{h,c})(\phi))=\TT_{h,c}(E_{\Bun,1/2,n}(\phi)).
\end{equation}
\subsubsection{Eigen-functions and cuspidal functions: the idea}

Define
\begin{definition}
\begin{enumerate}
\item A function $\phi\in \calS(\Bun_n)_\chi$ is cuspidal if the span of $\{T_{h,n}(\phi)\}, h\in \calH_{c}(G) , c\in \mcC _n$ is finite dimensional.
\item $ \calS_{cusp}(\Bun_n) \subset \calS(\Bun_n) $ is the subspace of cuspidal functions. \end{enumerate}

\end{definition}

\begin{conjecture}\label{cusp}$\calS_{cusp}(\Bun_n) _\chi $ is finite-dimensional for any $n$ and $dim( \calS_{cusp}(\Bun_n)_\chi)\sim q^{n\, dim \rm {Bun}}$ for $q>>1$.
\end{conjecture}

The equation (\ref{eisen-n}) implies that if $\phi$ is an eigen-function for all the operators $\gam_{c,n}(\TT_{h,c})$ (for
all $c$) then $E_{\Bun,1/2,n}(\phi)$ is an eigen-function for $\TT_{h,c}$. In the case $n=1$ such eigen-functions $\phi$ are described by the Langlands conjectures. We would like to explain why there are many such eigen-functions for $n=2$ and $G=\GL_2$ (but we still don't know how to describe the corresponding eigen-values).
Our strategy will be as follows: we are going to define an analog of cuspidal functions. This notion will be invariant under the operators $\gam_{c,n}(\TT_{h,c})$ and the space of cuspidal functions will be finite-dimensional. It will be easy to show that the operators $\gam_{c,n}(\TT_{h,c})$ are diagonalisable in the space of cuspidal functions -- in this way we'll construct a lot of eigen-functions.

\subsection{The case of $G=\GL_2$}In this subsection we outline a proof Conjecture \ref{cusp} in the case when $G=\GL_2,n=2$. We start with a different (but equivalent) definition of cuspidality.

\subsubsection{Definition of cuspidal functions}
We would like to define an analog of constant term operator in our setting.
Recall that the usual constant term operator (for $n=1$) is defined as follows. Let $P$ be a parabolic subgroup of $G$; it has a natural homomorphism to $M$ -- the Levi factor. Consider the diagram
\begin{equation}\label{ct-finite}
\CD
\Bun_P(k) @>{p}>> \Bun_G(k) \\
@V{q}VV  \\
\Bun_M(k).
\endCD
\end{equation}
Then the constant term $c_{G,P}$ is equal to $q_!p^*$. And a function $\phi$ on $\Bun_G(\mcC_k)$ is cuspidal if $c_{G,P}(\phi)=0$ for all parabolic subgroups $P$ of $G$; let $\calS_{cusp}(\Bun_G(k))$ be the space of all cuspidal functions. Then the following is well-known and is easy to prove:

1) $\calS_{cusp}(\Bun_G(k))$ is invariant under the Hecke operators.

2) $\calS_{cusp}(\Bun_G(k))$ consists of functions with finite support.

3) $\dim \calS_{cusp}(\Bun_G(k))<\infty$.

We would like to define an analogous space with the above properties for $\Bun_G(\mcC_n)$.
For this we need the analog of the constant term operator $\bfc_{G,P}$ for all $n$ which we would like to denote
by $\bfc_{G,P}^{(n)}$.
The definition is not completely straightforward; we shall only discuss the case $G=\GL_2$ and $n=2$.
So, we shall now assume that $G=\GL_2$ and we shall again just write $\Bun$ instead of $\Bun_{\GL_2}$. Also in this case the only proper parabolic subgroup up to conjugacy is the Borel subgroup; we shall also denote the corresponding constant term operator (that we are going to define) simply by $\bfc^{(2)}$.

In this case $T=\mathbb G_m\times \mathbb G_m$, so $\Bun_T(\mcC_2)=\Pic(\mcC_2)\times \Pic(\mcC_2)$, so it would be natural to expect that our constant term operator $\bfc^{(2)}$ maps functions on $\Bun(\mcC_2)$ to functions on $\Pic(\mcC_2)\times \Pic(\mcC_2)$. However, we do not know how to define such an operator if we want it to commute with the Hecke operators in some reasonable sense.
Instead, let us do the following.
Consider the semi-group $\Pic_2'$ (which contains the Picard group $\Pic_2$ of $\mcC_2$.
By definition, $\Pic_2'$  consists of coherent sheaves $\mcM$ on $\mcC$ such that $t \mcM \neq 0$ and
 there exists an imbedding $\mcM \ho \mcL$ where $\mcL$ is a line bundle on $\mcC$. The tensor product defines the semi-group structure on $\Pic_2'$.

\begin{example} Let $\mcC = \Spec (A_2[x]), J\subset A_2[x]$ be the maximal ideal generated by $(x,t)$ and $\mcJ$ the corresponding sheaf on $\mcC$. Then $ \mcJ \otimes \mcJ = \mcI$ where $I\subset xA_2[x]$ is generated by $(x^2,tx)$.
\end{example}

We would like now to define an analog of the diagram (\ref{ct-finite}). Namely, we consider the diagram
\begin{equation}\label{ct-2}
\CD
\Bun_B'(\mcC_2) @>{p_2}>> \Bun_2 \\
@V{q_2}VV  \\
\Pic'_2\times \Pic'_2.
\endCD
\end{equation}
where $\Bun_B'(\mcC_2)$ consists of all short exact sequences
$$
0\to\mcL_1\to \mcF\to \mcL_2\to 0
$$
where $\mcF \in  \Bun_2, \mcL_1,\mcL_2 \in \Pic_2 '$. It is easy to see that in this case we have $\mcL_i=\underline{\Hom}(\mcL_j,\det(\mcF))$ for $i,j=1,2$ and $i\neq j$. So if we fix $\det(\mcF)$ and one of the bundles
$\mcL_1$ of $\mcL_2$, this determines the isomorphism class of the other.
\begin{definition}
\begin{enumerate}
\item
Set $\bfc^{(2)}=(q_2)_!p_2^*$.
\item
We denote by $\calS_{cusp}'(\Bun_2)$ the space of all $\phi\in \calS(\Bun_2)$ such that $\bfc^{(2)}(\phi)=0$.
\end{enumerate}
\end{definition}

We now claim the following:
\begin{theorem}
\begin{enumerate}
\item
The space $\calS_{cusp}'(\Bun_2)$ is invariant under the Hecke operators - i.e. it is preserved by all the operators of the form $\gam_{c,n}(\TT_h)$.
\item
$\calS_{cusp}'(\Bun_2)=\calS_{cusp}(\Bun_2)$.
\item
Fix a character $\chi:\Pic_2\to \mathbb C^*$. Let $\calS_{cusp,\chi}(\Bun_2)$ denote the space of cuspidal functions
on $\Bun_2$ s.t. $\phi(\mcF\otimes \mcL)=\phi(\mcF)\cdot \chi(\mcL)$ where $\mcF\in \Bun_2, \mcL\in \Pic_2$ and such that the support of $\phi$ modulo $\Pic_2$ is finite (the last condition probably follows from cuspidality).
Then $\dim \calS_{cusp,\chi}(\Bun_2)<\infty$.
\end{enumerate}
\end{theorem}
The proof will be discussed in another publication. Let us just comment on the first assertion.
In this case the algebra $\calH_c(G)$ is generated by three elements - corresponding to coweights $\lam=(1,0), \mu=(1,1)$ and $-\mu=(-1,-1)$ of $\GL(2)$. The statement is essentially trivial for the Hecke operators corresponding to
$\mu$ and $-\mu$.
In case of the operator $\gam_{c,2}(\TT_{\lam})$ we can costruct an explicit operator $B_c$ acting from the space of (all) functions on $\Pic_2'\times \Pic_2'$ such that $\bfc^{(2)}\circ \gam_{c,2}(\TT_{\lam,c})=B_c\circ \bfc^{(2)}$. This will imply that the kernel of $\bfc^{(2)}$ is invariant under the operators $\gam_{c,2}(\TT_{\lam,c})$.

\subsection{Main question}
It would be extremely interesting to  relate the spectrum of the algebra $\calH _2(Gl_2)$ acting on  $\mcS _{cusp ,\chi}(\rm {Bun_2})$
to some kind of Galois data (involving the dual group $G^{\vee}$). At the moment we do not know how to do it
in the non-archimedian case. In the archimedian case we formulated a precise conjecture of this sort in \cite{efk1} and \cite{efk2}.
\section{The case $F=\mC$}\label{complex}

\subsection{From Hecke operators to differential operators: the idea}
In this section we specialize to the case $F=\CC$. In this case in addition to Hecke operators one can introduce another player: the algebra of twisted (polynomial) differential operators on $\Bun$, which will roughly speaking act on the same space as the Hecke operators and the two actions will commute. This will allow us to formulate a variant of Langlands conjectures in this case. More precisely, we are going to relate the Hecke eigen-values to some particular $G^{\vee}$-local systems on $\mcC$ -- opers with real monodromy .
Let us begin by recalling basic information about opers and differential operators on $\Bun$.

\subsection{Opers}
 For a principal $G^\vee$-bundle $\mcG$ on $\mcC$ we denote by $\Fl_\mcG$ the associated $\Fl$- bundle on $\mcC$ where $Fl$ is the variety of Borel subgroups  of $G^\vee$.

\begin{definition}\begin{enumerate}
\item
An $G^\vee $-oper on $\mcC$ is triple $(\mcG ,\nabla ,s)$ where
$\mcG$ is a principal $G^\vee$-bundle on $\mcC, \nabla$ is a connection on $\mcG$ and $s$

is a section of $\Fl _\mcG$ satisfying an analog of the   Griffiths type condition with  respect to $\nabla$ (see \cite{bd1}). We denote by $\Oper_{G^\vee}(\mcC)$ the variety of opers.
\footnote{If $G^\vee$ is adjoint then the moduli stack of opers is in fact an algebraic variety (which is isomorphic to an affine space of dimension $\rk(G)$). If $G^\vee$ is not adjoint then formally one needs to consider the coarse moduli space here, since the center $Z$ of $G^\vee$ is equal to the group of automorphisms of every oper. We shall ignore this subtlety for the rest of this section.}
\item For an oper  $o= (\mcF ,\nabla ,s)$ we denote by $\zeta _o: \pi _1(C)\to G^\vee (\mC)$ the morphism defined by the connection $\nabla$. We denote by $\Oper _G^\vee (\mcC)^\mR \subset \Oper _G^\vee (\mcC) $ the subset of opers $o$ such that the homomorphisms $ \zeta _o $ and $\ov\zeta _o $ are conjugate, where $^{-}  : G^\vee (\CC) \to G^\vee (\CC)$ is the complex conjugation corresponding to a choice of a split real form of $G^\vee$.
\end{enumerate}
\end{definition}
Let us make several comments. First, it is known (cf. \cite{bd}) that given just a pair $(\mcG,\nabla)$, the $B$ structure $s$ is unique if it exists. Thus $\Oper_{G^\vee}(\mcC)$ is actually a closed subset of the moduli stack of $G^\vee$-bundles with a connection (in other words, for such a local system to be an oper is a property rather than a structure). Second, let us comment on the reality condition in (2). Obviously, one way to guarantee this condition is to require that the monodromy representation of $\pi_1(\mcC)$ corresponding to $(\mcG,\nabla)$ is conjugate to a homomorphism going into $G ^\vee (\mR)$ for a real split form of $G ^\vee $. We expect that the converse is also true, and this is proved for $G ^\vee =SL(2)$ in \cite{efk2} (Remark 1.8), but we don't know how to prove this in general. However, it is not hard to see (cf. again \cite{efk2}) that up to conjugation the image of the monodromy homomorphism $\pi_1(\mcC)\to G ^\vee $ lies in some inner form of the split real form of $G ^\vee $. When we are in the setup of subsection \ref{parabolic} and $|D|\geq 1$ it is also shown in {\em loc. cit.} that the monodromy is lies in the split real form of $G ^\vee $.
%--------------------------------------------------------------------------------------
\subsection{Opers and differential operators}
Let $\mcD$ be the algebra of global sections of the sheaf $\DD_{1/2}(\Bun)$ of regular differential operators on $\bo _{\Bun}^{1/2}$. We denote by $\tau:\calD\to \calD$ the involution on $\calD$ induced by the Cartan involution of $G $.

The following statement is one of the main results of \cite{bd} (a local version of this result appears in \cite{ff}).
\begin{theorem}\label{b}
\begin{enumerate}
\item The algebra $\mcD$ is commutative.
\item $\Spec(\mcD )= \Oper _{ G^{\vee}}(\mcC) $.
\item Let $o\in \Oper _{ G^{\vee}}(\mcC)$ and let $\chi_o:\calD\to \CC$ be the corresponding homomorphism. Let also $I_o\subset \DD_{1/2}(\Bun)$ be the sheaf of ideals of $\DD_{1/2}(\Bun)$ generated elements of the form $d-\chi_o(d)$ where $d\in \calD$. Then the $\DD_{1/2}(\Bun)$-module $M_o:=\DD_{1/2}(\Bun)/I_o$ is $\calO_{\Bun}$-coherent when restricted to $\Bun_{vs}$.
\begin{footnote}{Part (3) of this Theorem explains the reason for our belive in Conjecture  \ref{spaces}(2)}\end{footnote}

\end{enumerate}
\end{theorem}

\subsection{Differential operators and Hecke operators}
Recall that we denote by $C^{\infty}_{1/2}(\Bun_{vs}) $ be the space of smooth
$1/2$-forms on  $ \Bun_{vs}$. The algebra $\calA := \calD \otimes \ov {\calD}$ acts naturally on   $ C^{\infty}_{1/2}(\Bun_{vs}) $. We denote by $\hat\tau$ the involution on $\calA$ such that $\hat \tau (d_1\otimes \bar d_2)= d_2 ^\tau \otimes \ov {d_1^\tau}$ and define
$\calA ^\mR \subset \calA $ as the subalgebra of $\hat \tau$-fixed points.

We would like to claim that the action of the algebra $\calA$ on $1/2$-forms commutes with the action of the Hecke operators. Here we must be careful, as a priori it is not clear how to construct one vector space on which both algebras will act. For this we need to formulate one more definition.

Let us define a space $\Sch(\Bun)$ -- "the Schwartz space of $\Bun$". Namely, we set
\begin{equation}\label{sch}
\Sch(\Bun)=\{ \phi\in C^{\infty}_{1/2}(\Bun_{vs})|\ a(\phi)\in L^2(\Bun)\ \text{for any $a\in \calA$}\}
\end{equation}
For $a\in \calA$ we denote by $\hat a$ the induced endomorphism of $\Sch(\Bun)$.

Note that by definition $\Sch(\Bun)\subset L^2(\Bun)$ and also $\calS(\Bun_{vs})\subset \Sch(\Bun)$. The reader might ask why we start with $C^{\infty}$-forms on $\Bun_{vs}$ rather than on $\Bun_{st}$. The reason is that below we want to study eigen-vectors of $\calA$ on $\Sch(\Bun)$ and it follows easily from Theorem \ref{b}(3) that any such eigen-vector is automatically smooth on $\Bun_{vs}$ (but there is no reason for it to be smooth on $\Bun_{st}$).

\begin{conjecture}\label{4}

\begin{enumerate}
\item Any $a\in \calA ^\mR$ extends to an (unbounded) self-adjoint operator on $L^2(\Bun)$.
\item
The space $\Sch(\Bun)$ is stable under the action of all Hecke operators.
\item
$\Sch(\Bun)=\iota_{1/2}(\calS_{1/2}(\Bun)$.
\footnote{Note that if we assume the validity of Conjecture \ref{spaces} for $F=\CC$, then (3) implies (2).}
\item
The action of $\calA$ on $\Sch(\Bun)$ commutes with the action of Hecke operators.
\item
There exists a dense (in the $L^2$-sense) subspace $\Sch(\Bun)_0$ of $\Sch(\Bun)$ which is stable under  $\calA$ and
the Hecke operators and such that $\Sch(\Bun)_0$ is a direct sum of 1-dimensional eigen-spaces for $\calA$  (in other words, the space $\Sch(\Bun)_0$ is locally finite-dimensional for $\calA$ and every generalized eigen-value has multiplicity 1).
\end{enumerate}
\end{conjecture}

Let us comment on the multiplicity 1 statement. A 1/2-form is actually an eigen-vector if it satisfies certain system of linear differential equations. Locally on $\Bun_{vs}(\CC)$ the space of solutions is finite-dimensional but certainly not one-dimensional (this has to do with the fact the $\DD$-module $M_o$ has high rank on $\Bun_{vs}$; for example, for $\SL_2$ this rank is $2^{3g-3}$). However, globally most of these solutions become multi-valued, so the multiplicity one conjecture says that only one-dimensional space of solutions is single-valued globally. This in fact would follow if we knew that the $\DD$-module $M_o$ was irreducible and had regular singularities.  For $G=\PGL_2$ this can be deduced from \cite{Gai} (and probably similar analysis can be carried over for $\PGL_n$).

Conjecture \ref{4} implies that $L^2(\Bun)$ is a (completed) direct sum of eigen-spaces for
$\calA$ and eigen-values have multiplicity 1. A priori any such eigen-value is given by a pair of opers $(o,o')$, but the part $(1)$ of Conjecture \ref{4} implies that $o'=\ov {o^\tau}$, so we are supposed to attach an eigen-space to a single oper $o$. It is also not difficult to see that $o\in \Oper_{G^{\vee}}(\mcC)^{\mR}$. We denote the corresponding eigen-space by
$L^2(\Bun)_o$. Note that Conjecture \ref{4} implies that $L^2(\Bun)_o\subset \Sch(\Bun)$.

\begin{conjecture}\label{real-oper}
We have $L^2(\Bun)_o\neq 0$ if and only if $o\in \Oper_{G^{\vee}}(\mcC)^{\mR}$.
\end{conjecture}

\begin{remark}
As was remarked above, the "only if" direction is easy. What is hard is to prove existence of eigen-vectors for $\calA$ which lie in $L^2$.
\end{remark}

Note  that Conjectures \ref{compact}, \ref{4} and \ref{real-oper} together imply the following
\begin{corollary}\label{corol}
Let $\mcW$ denote the set of Hecke eigen-values on $L^2(\Bun)$. Then
there exists a surjective map $\eta :\Oper _{ G^{\vee}}(\mcC) ^\mR \to  \mcW $ such that for any $c\in \mcC$ and any
$h\in \calH(G,\CC)$ the operator
$\TT_{h,c}$ acts on $L^2(\Bun)_0$ by  $\eta (o)(h)$.
\end{corollary}

Let us comment on the connection between Corollary \ref{corol} and Conjecture \ref{compact}. We actually expect the map $\eta$ to be finite-to-one (and in many cases it should be an isomoprhism), so Conjecture \ref{compact} should imply that
$\Oper _{ G^{\vee}}(\mcC) ^\mR$ should be a discrete subset of $\Oper _{ G^{\vee}}(\mcC)$. This assertion is not obvious, and at the moment we don't know how to prove it in general, but let us note that for $G^{\vee}=\PGL_2$ it was proven by G.~Faltings in \cite{fa-real}.

\subsection{Eigenvalues of Hecke operators} We conclude this section by describing a conjectural formula for the map $\eta$.
(the contents of this subsection are described in more detail in \cite{efk2}).
More precisely, we are going to do the following. We would like to understand the scalar by which the operator $\TT_{\lam,c}$ acts in $L^2(\Bun)_o$. We can actually regard $c$ as a variable here. In view of Remark \ref{hecke-omega} this eigen-value is in fact a section $\Phi_{\lam,o}$ of $|\ome_{\mcC,c}|^{-\langle \lam,\rho^{\vee}\rangle}$ (recall that $\rho^{\vee}$ denotes the half-sum of positive coroots of $G$).

For $\bl \in \bL ^+$, let $V_\bl$ be the corresponding
irreducible finite-dimensional representation of $G^{\vee}$. Choose an $o=(\mcF,\nabla,s)\in \Oper _{ G^{\vee}}(\mcC)$. Moreover, the Griffiths transversality condition impies that the $T^{\vee}$-bundle induced from the $B^{\vee}$-structure $s$ by means of the homomorphism $B^{\vee}\to {T^{\vee}}$ is induced from $\ome_{\mcC}$ by means of the cocharacter $\rho^{\vee}:{\mathbb G_m}\to {T^{\vee}}$
\footnote{Strictly speaking, this makes sense only if $G^{\vee}$ is simply connected since $\rho^{\vee}$ is a well-defined cocharacter of $T^{\vee}$ only in that case. For general $G$ the corresponding $T^{\vee}$-bundle is induced from $\ome_{\mcC}^{1/2}$ by the character $2\rho^{\vee}$ for some choice of $\ome_{\mcC}^{1/2}$. To simplify the notation we are going to write the answer in the case when $G^{\vee}$ is simply connected -- the generalization to any $G$ is straightforward}.
Therefore if
we denote by $(\V_{o,\bl},\nabla_{o,\bl})$ the vector bundle on $\mcC$ associated to $\mcF$ via the representation $V_{\lam}$ (with the corresponding flat connection), then $s$ defines
an embedding
  $$
  \bo _C^{\langle \bl,\rho^{\vee} \rangle} \hookrightarrow
  \V_{o,\bl}
  $$
  and hence a morphism
  $$
 \mcO _{\mcC} \hookrightarrow \bo _C^{-\langle \bl,\rho^{\vee} \rangle} \otimes
  \V_{o,\bl}.
  $$
  We let $\sigma_{\lam}$ be the image of 1 under this morphism.

Let now $o\in \Oper _{ G^{\vee}}(C )^\mR$. Then we have
isomorphism of $\V_{o,\lam}$ and $\ol{\V}_{o,\lam}$ as flat $C^\infty$-bundles
(and this isomorphism is canonical up to the action of the center of $G^{\vee}$).
Since $V^*_\la \simeq V_{-w_0(\la)}$ we get a pairing $(\cdot,\cdot)_{\lam}$ between $C^{\infty}$-sections of
$\V_{o,\lam}$ and of  $\ol{\V}_{o,\lam}$ .
Since $\langle
  -w_0(\bl),\rho^{\vee} \rangle = \langle \bl,\rho^{\vee} \rangle$, we can regard
$\ol{\sigma}_{-w_0(\bl)}$ as a section of  $\ol{\bo}_{\mcC}^{-\langle
    \bl,\rho^{\vee} \rangle} \otimes \ol{\V}^*_\bl)$. Since\linebreak
    $\ome_{\mcC}^{-\langle \bl,\rho^{\vee} \rangle}\otimes \ol{\ome}_{\mcC}^{-\langle \bl,\rho^{\vee} \rangle}=|\ome_{\mcC}|^{-\langle \bl,\rho^{\vee} \rangle}$ we can formulate the following Conjecture (cf. \cite{efk2}):
\begin{conjecture}
$$
\Phi_{\lam,o}=(\sigma_{\lam},\ol{\sigma}_{-w_0(\lam)})_{\lam}\in C^{\infty}(\mcC,|\ome_{\mcC}|^{-\langle \bl,\rho^{\vee} \rangle}).
$$
\end{conjecture}
\subsection{Parabolic bundles: results} All the conjectures of this Section can be easily generalized to the setup of subsection \ref{parabolic}. In the case when $\mcC$ is $\PP^1$ and the cardinality of the divisor $D$ is 3, 4 or 5 they are proven in \cite{efk3} (and most of them are proven in {\em loc. cit.} even for $|D|>5$).

\section{The case $F=\Bbb R$}\label{real}

In this section we would like to describe the conjectural picture of the analytic Langlands correspondence in the case $F=\Bbb R$. This picture has been developed by P. Etingof, E. Frenkel, D. Gaiotto, D. Kazhdan and E. Witten, and is discussed in \cite{gw}, Section 6.

\medskip
\noindent
{\bf Warning.} Some of the letters used in the previous section (such as $\sigma$ or $\tau$) will have a different meaning in this section.

\subsection{Real groups, $L$-groups and all that}
Let $G$ be a connected complex semi-simple group. Recall that a {\it real structure} on $G$ is defined by an antiholomorphic involution $\sigma: G\to G$. The corresponding {\it group of real points} is $G^\sigma$ (it may be disconnected). The inner class of $\sigma$ gives rise to a based root datum involution $s=s_\sigma$ for $G$ which is also one for $G^\vee$. If $G$ is semisimple, this is just a Dynkin diagram automorphism.

Recall \cite{ABV} that to $G,s$ we may attach the {\it Langlands L-group} ${}^LG={}^LG_s$, the semidirect product of $\Bbb Z/2={\rm Gal}(\Bbb C/\Bbb R)$ by $G^\vee$, with the action of $\Bbb Z/2$ defined by $\gamma\circ s$, where $\gamma$ is the Cartan involution.
%--------------------------------------------------------------------------------------------------------
\subsection{$L$-systems}
Let $\mcC$ be a compact complex Riemann surface of genus $g\ge 2$. Let $\tau: \mcC\to \mcC$
be an antiholomorphic involution.
Given a holomorphic principal $G$-bundle $\calE$ on $\mcC$, we can define
the antiholomorphic bundle $\tau(\calE)$, hence a holomorphic bundle
$\sigma \tau(\calE)$. Let us say that $\calE$ is {\it real} under $\sigma$
if there exists an isomorphism  $A: \calE\to \sigma\tau(\calE)$ such that
\begin{equation}\label{condi}
\sigma\tau(A)\circ A=1.
\end{equation}
This isomorphism $A$ is unique if exists and \eqref{condi} is automatic if ${\rm Aut}(\calE)=1$, which happens generically for stable bundles if $G$ is adjoint.
In this case
$gA: \calE\to g\sigma\tau(\calE)$ has the same property for $\sigma'=g\sigma$, where $g\in  G$ and $g\sigma(g)=1$. Thus the moduli space
of such stable bundles depends only on $s$ (\cite{BGH}, Proposition 3.8). We will denote it by $\Bun_{G,s}$.

Consider the simplest case when $\tau$ has no fixed points, i.e., $\mcC(\Bbb R)=\emptyset$. Let $\zeta$ be a local system on the non-orientable surface $\mcC/\tau$ with structure group ${}^LG$. Let us say that $\zeta$ is an {\it L-system} if it attaches to every orientation-reversing path in $\mcC/\tau$ a conjugacy class in ${}^LG$ that maps to the nontrivial element in $\Bbb Z/2$. The following conjecture is formulated in \cite{gw}, Section 6 (in the case of the compact inner class).

\begin{conjecture} The spectrum of Hecke operators on $L^2(\Bun)$
is parametrized by $L$-systems on $\mcC/\tau$ with values in ${}^LG={}^LG_s$
whose pullback to $\mcC$ have a structure of a $G$-oper.
\end{conjecture}

\begin{example} Let $s=\gamma$. Then ${}^LG=\Bbb Z/2\times G^\vee$, so
an L-system is the same thing as a $G^\vee$-local system on $\mcC/\tau$.
So in this case the condition on the $G^\vee$-local system on $\mcC$ to occur in the spectrum
is (conjecturally) that it extends to the 3-manifold $M:=(\mcC\times [-1,1])/(\tau, -{\rm Id})$ whose boundary is $\mcC$ (and this extension is a part of the data).

Namely, in this case the spectral local systems are $\zeta$ which are isomorphic to $\zeta^\tau$ and such that $\zeta$ is an oper (hence also an anti-oper), so $\zeta$ is a real oper
``with real coefficients". But among these we should only choose those that extend
to $\mcC/\tau$ (and then the multiplicity of eigenvalue may be related to the number of such extensions). This agrees with the picture \cite{gw}, 6.2 coming from 4-dimensional supersymmetric gauge theory.\footnote{More precisely, as was explained to us by E. Witten, what comes from ordinary gauge theory is this picture for the compact inner class $s$. To obtain other inner classes, one needs to consider {\rm twisted gauge theory} where the twisting is by a Dynkin diagram automorphism of $G$. Namely, gauge fields in this theory are
invariant under complex conjugation $\tau$ up to such an automorphism.}
More precisely, recall that by a result of Beilinson and Drinfeld \cite{bd}, opers for adjoint groups have no nontrivial automorphisms. So for any connected semi-simple $G$ we get an obstruction for such $\zeta$ to extend to $\mcC/\tau$ which lies in $Z/Z^2=H^2(\Bbb Z/2,Z)$, where $Z$ is the center of $G^\vee$.\footnote{Indeed, $\pi_1(\mcC/\tau)$ is generated by
$\pi_1(\mcC)$ and an element $t$ such that $tbt^{-1}=\beta(b)$
for some automorphism $\beta$ of $\pi_1(\mcC)$, and
$t^2=c\in \pi_1(\mcC)$, so that $\beta^2(b)=cbc^{-1}$.
So given a representation $\zeta: \pi_1(\mcC)\to G^\vee$,
an L-system would be given by an assignment $\zeta(t)=T\in G^\vee$ such that
(1) $T^2=\zeta(c)$ and (2) $T\zeta(a)T^{-1}=\zeta(\beta(a))$.
If $\zeta\cong \rho\circ \beta$ then
$T$ satisfying (2) is unique
up to multiplying by $u\in Z$, and $T^2=\zeta(c)z$, $z\in Z$.
Moreover, if $T$ is replaced by $Tu$ then $z$ is replaced by $zu^2$, hence
the obstruction to satisfying (1) lies in $Z/Z^2$.}
Moreover, if this obstruction vanishes then the freedom for choosing the extension  is in a torsor over $H^1(\Bbb Z/2,Z)=Z_2$, the 2-torsion subgroup in $Z$.
\end{example}

\begin{example} Let $G=K\times K$ for some complex group $K$, and
$s$ be the permutation of components (the only real form  in this inner class is $K$
regarded as a real group). This is equivalent to the case $F=\Bbb C$ considered above (for $\mcC$ defined over $\Bbb R$). Then ${}^LG_s={}^LG_{\gamma\circ s}=\Bbb Z/2\ltimes (K^\vee\times K^\vee)$, where $\Bbb Z/2$ acts by permutation. So an $L$-system is a $K^\vee\times K^\vee$ local system on $\mcC$ of the form
$(\rho,\rho^\tau)$. Thus the spectrum is parametrized by $\zeta$ such that both $\zeta$ and $\zeta^\tau$ are opers, i.e., $\zeta$ is both an oper and an anti-oper, i.e. a real oper, which agrees with the conjecture for $F=\Bbb C$. (Note that in this case
$H^i(\Bbb Z/2,Z)=1$ so there is no obstructions or freedom for extensions).
\end{example}

\begin{remark} If $\mcC(\Bbb R)\ne \emptyset$, the story gets more complicated,
and we will not discuss the details here. Let us just indicate that, as explained in \cite{gw}, Section 6, to define the appropriate moduli space and the spectral problem on it, we need to fix a real form $G_i$ of $G$ in the inner class $s$ for each component (oval) $C_i$ of $\mcC(\Bbb R)$,
and the eigenvalues of Hecke operators are conjecturally parametrized by a certain kind of ``real" opers corresponding to this data, i.e., opers with real coefficients satisfying appropriate reality conditions on the monodromy of the corresponding $G^\vee$-connection.
Furthermore, in the tamely ramified case, when we also have a collection of marked points $D$ on $\mcC$ defined over $\Bbb R$, to define the most general version of our spectral problem, we need to fix a unitary representation $\pi_i$ of the real group $G_i$ for every marked point $c\in D$ on $C_i$ and a unitary representation of the complex group $G_\Bbb C$ for every pair of complex conjugate marked points $c,\overline c\in D$ not belonging to $\mcC(\Bbb R)$. For example, the case of parabolic structures corresponds to taking $s$ to be the split inner class, $G_i$ the split forms, and  $\pi_i$ the unitary principal series representations. In the genus zero case,
this was discussed in detail in \cite{efk3}, and it was shown that this problem leads to appearance of $T$-systems.
\end{remark}

\subsection{Connection to Gaudin model}
Recall that the {\it Gaudin model} for a simple complex Lie algebra $\grg$ is the problem of diagonalization of the {\it Gaudin hamiltonians}
$$
H_i:=\sum_{1\le j\le N,j\ne i}\frac{\Omega_{ij}}{z_i-z_j}
$$
on the space $(V_1\otimes...\otimes V_N)^\grg$, where $V_i$ are finite dimensional $\grg$-modules, $z_i\in \Bbb C$ are distinct points, $\Omega\in (S^2\grg)^\grg$ is the Casimir tensor dual to the Killing form, and $\Omega_{ij}$ denotes the action of $\Omega$ in the $i$-th and $j$-th factor. These operators commute, and if $\grg\ne \mathfrak{sl}_2$ then there are also higher Gaudin hamiltonians associated to the Feigin-Frenkel higher Sugawara central elements at the critical level (see \cite{FFR}), which commute with each other and with $H_i$, and the problem is to simultaneously diagonalize all these operators.

It turns out that this problem (for real $z_i$) is a special case of the spectral problem considered in this paper, in the case $F=\Bbb R$. Namely, let us take $\mcC=\Bbb P^1$ with the usual real structure and fix the compact inner class $s$ of the complex simply connected group $G$ with ${\rm Lie}(G)=\grg$. As explained in the previous remark, on the real locus $\Bbb P^1(\Bbb R)$ we are supposed to fix a real form of $G$ in this inner class, and we fix
the compact form $G_c$. Further, consider marked points $z_1,...,z_N$ on the real locus
(the tamely ramified case). Then we are supposed to fix a unitary representation
of $G_c$ at every $z_i$, and we take it to be $V_i$. \linebreak  Then the Hilbert space
 of the analytic Langlands theory is \linebreak  $\mathcal H=(V_1\otimes...\otimes V_n)^{G_c}$ (so in this case it is finite dimensional), and the quantum Hitchin system is comprised exactly by the Gaudin hamiltonians (including the higher ones), cf. \cite{FFR}.

As is explained in \cite{F1,F2}, the Bethe ansatz method shows that the eigenvectors of the Gaudin hamiltonians are labeled by monodromy-free $G^\vee$-opers on $\Bbb P^1$ with first order poles at $z_i$ and residues in the conjugacy class of $-\lambda_i-\rho$, where $\lambda_i$ is the highest weight of $V_i$. These are exactly the ``real opers" for this situation. Thus the results of \cite{F1,F2} may be considered as a finite-dimensional instance of the tamely ramified analytic Langlands correspondence for genus zero and $F=\Bbb R$.

\end{document}